\documentclass[12pt]{article}
\usepackage{amsmath,amssymb}
\newcommand{\lan}{\langle}
\newcommand{\ran}{\rangle}
\newcommand{\lt}{\left}
\newcommand{\rt}{\right}
\newcommand{\pf}{\mbox{\bf Proof}}
\newcommand{\ra}{\rightarrow}
\newcommand{\qed}{\hfill{\rule {.1in}{.1in}}}

\newcommand{\ol}{\overline}

\newcommand{\vp}{\varphi}
\newtheorem{thm}{Theorem}[section]
\newtheorem{lma}[thm]{Lemma}
\newtheorem{prop}[thm]{Proposition}
\newtheorem{cor}[thm]{Corollary}
\title{$\lambda$-Toeplitz operators with analytic
symbols}
\author{Chih Hao Chen, Po Han Chen, Mark C. Ho\thanks{The research is part of student's thesis, under the supervision of Mark C. Ho} and Meng Syun Syu\\
\scriptsize{Department of Applied Mathematics, National Sun Yat-Sen
University, Kaohsiung, Taiwan 80424}}
\date{\scriptsize{{\bf AMS classification}: 47B33, 47B35}}
\begin{document}
\maketitle
\begin{abstract}
Let $\lambda$ be a complex number in the closed unit disc $\ol{\Bbb
D}$, and $\cal H$ be a separable Hilbert space with the orthonormal
basis, say, ${\cal E}=\{e_n:n=0,1,2,\cdots\}$. A bounded operator
$T$ on $\cal H$ is called a {\em $\lambda$-Toeplitz operator} if
$\lan Te_{m+1},e_{n+1}\ran=\lambda\lan Te_m,e_n\ran$ (where
$\lan\cdot,\cdot\ran$ is the inner product on $\cal H$). The
subject arises naturally from a special case of the operator
equation
\[
S^*AS=\lambda A+B,\ \mbox{where $S$ is a shift on $\cal H$},
\]
which plays an essential role in finding bounded matrix $(a_{ij})$
on $l^2(\Bbb Z)$ that solves the system of equations
\[\lt\{\begin{array}{lcc}
a_{2i,2j}&=&p_{ij}+aa_{ij}\vspace{.05in}\\
a_{2i,2j-1}&=&q_{ij}+ba_{ij}\vspace{.05in}\\
a_{2i-1,2j}&=&v_{ij}+ca_{ij}\vspace{.05in}\\
a_{2i-1,2j-1}&=&w_{ij}+da_{ij}
\end{array}\rt.
\]
for all $i,j\in\Bbb Z$, where $(p_{ij})$, $(q_{ij})$, $(v_{ij})$,
$(w_{ij})$ are bounded matrices on $l^2(\Bbb Z)$ and $a,b,c,d\in\Bbb
C$. In this paper, we study the essential spectra for $\lambda$-Toeplitz operators when $|\lambda|=1$, and we will use the results to determine the spectra of certain weighted composition operators on Hardy spaces.
\end{abstract}

\section{Introduction}

Let $\cal H$ be a separable Hilbert space with an orthonormal basis,
say, ${\cal E}=\{e_n:n=0,1,2,\cdots\}$. Given $\lambda\in\ol{\Bbb
D}=\{z\in\Bbb C:|z|\leq1\}$, a bounded operator $T$ is called a {\em
$\lambda$-Toeplitz operator} if $\lan
Te_{m+1},e_{n+1}\ran=\lambda\lan Te_m,e_n\ran$ (where
$\lan\cdot,\cdot\ran$ is the inner product on $\cal H$). In terms of
the basis $\cal E$, it is easy to see that the matrix representation
of $T$ is given by
\[\lt(\begin{array}{cccccc}
a_0&a_{-1}&a_{-2}&a_{-3}&a_{-4}&\cdots\\
a_1&\lambda a_0&\lambda a_{-1}&\lambda a_{-2}&\lambda a_{-3}&\ddots\\
a_2&\lambda a_1&\lambda^2a_0&\lambda^2a_{-1}&\lambda^2a_{-2}&\ddots\\
a_3&\lambda a_2&\lambda^2a_1&\lambda^3a_0&\lambda^3a_{-1}&\ddots\\
a_4&\lambda a_3&\lambda^2a_2&\lambda^3a_1&\lambda^4a_0&\ddots\\
\vdots&\ddots&\ddots&\ddots&\ddots&\ddots
\end{array}\rt)
\]
for some double sequence $\{a_n:n\in\Bbb Z\}$, and the boundedness
of $T$ clearly implies that $\sum|a_n|^2<\infty$. Therefore, it is
natural to introduce the notation
\[
T=T_{\lambda,\vp},
\]
where $\vp\sim\sum_{-\infty}^\infty a_ne^{in\theta}$ belongs to
$L^2=L^2(\Bbb T)$, the Hilbert space of square integrable functions
on the unit circle $\Bbb T$, with inner product
\[
\lan f,g\ran={1\over2\pi}\int_0^{2\pi}f{\ol g}d\theta,
\]
and consider $T_{\lambda,\vp}$ as an operator acting on $H^2$, i.e., the {\em Hardy space} defined by
\[
\lt\{f\in L^2:\int_0^{2\pi}f(e^{i\theta})e^{in\theta}d\theta=0,\ n<0
\rt\}
\]
 by identifying $\cal H$ with $H^2$ and $e_n$ identified with the function $e^{in\theta}$, $n\geq0$. Note that $H^2$ can be considered as the subspace of ``analytic functions" in $L^2$ since it consists of elements in $L^2$ such that the negative terms of their Fourier coefficients are zero. Also note that when $\lambda=1$ and $\vp\in L^\infty=L^\infty(\Bbb T)$, the matrix of $T_{1,\vp}$ is the matrix of the bounded Toeplitz operator $T_\vp$ on $H^2$. For the readers who are not familiar with the operator theory on $H^2$, the {\em Toeplitz operator $T_\vp$ with symbol} $\vp\in L^\infty$ is the operator defined by $T_\vp f=P(\vp f)$, $f\in H^2$, where $P$ is the projection from $L^2$ on to $H^2$. Here we refer the reader to \cite{RD:a} and \cite{KH:a}, both of which are  excellent sources of information on the theory of Hardy spaces and Toeplitz operators.

Our interests in $T_{\lambda,\vp}$ originated from the consideration
of the following action on ${\cal B}(\cal H)$: Let $S$ be the
unilateral shift, i.e., $Se_n=e_{n+1}$, $n=0,1,2,\cdots$, and define
the mapping on ${\cal B}(\cal H)$:
\[
\phi(A)=S^*AS,\ A\in{\cal B}(\cal H).
\]
Then it is not difficult to check, from the definition, that
$\phi(T_{\lambda,\vp})=\lambda T_{\lambda,\vp}$. Hence, the
$\lambda$-Toeplitz operators are precisely the ``eigenvectors" for
$\phi$ associated with $\lambda$, and notice that the Toeplitz operators are
just the special cases associated with $\lambda=1$. Moreover, since
$\|\phi\|\leq1$ (in fact, the spectrum of $\phi$ is the closed unit
disc), we accordingly restrict our attention to the case
$|\lambda|\leq1$. The motivation behind the action $\phi$ which
prompted our interests in this subject came from the fact that
this type of actions induced by shifts on Hilbert spaces is playing an important role in the study of the bounded
matrix $(a_{ij})$ on $l^2(\Bbb Z)$ (with respect to the canonical basis) which solves the system of equations
\[\lt\{\begin{array}{lcc}
a_{2i,2j}&=&p_{ij}+aa_{ij}\vspace{.05in}\\
a_{2i,2j-1}&=&q_{ij}+ba_{ij}\vspace{.05in}\\
a_{2i-1,2j}&=&v_{ij}+ca_{ij}\vspace{.05in}\\
a_{2i-1,2j-1}&=&w_{ij}+da_{ij}
\end{array}\rt.\eqno(\star)
\]
for all $i,j\in\Bbb Z$, where $(p_{ij})$, $(q_{ij})$, $(v_{ij})$,
$(w_{ij})$ are bounded matrices on $l^2(\Bbb Z)$ and $a,b,c,d\in\Bbb
C$ (In the analysis of this system, however, the shift involved has infinite multiplicity. See \cite{RR:a}). For
the details we refer the readers to \cite{Ho:f} and \cite{Ho:g}.

There is a big overlap between $\lambda$-Toeplitz operators and the
so-called {\em Toeplitz-composition operators}, i.e., operators
which can be expressed as the products of Toeplitz operators and {\em composition
operators} (See, for examples, \cite{MJ:a} and \cite{KMM:a}). To be specific, let
$|\lambda|\leq1$ and consider the operator $U_\lambda$ defined by
$U_\lambda e_n={\lambda}^ne_n$, $n=0,1,2,\cdots$. Then $U_\lambda$
can also be considered as the composition operator
\[
f\sim\sum_{n=0}^\infty a_ne^{i\theta}\ra f\circ\tau\sim\sum_{n=0}^\infty a_n\lambda^ne^{i\theta}
 \]
 on $H^2$ induced by the map $\tau(z)=\lambda z$, $z\in\Bbb D$ (or, with symbol $\tau$). Notice here we write the image of $f$ as $f\circ\tau$ because it is convenient to do so (Analytic maps from $\Bbb D$ into $\Bbb D$ always induce bounded composition operators on $H^2$. For the theory of composition operators, we refer the reader to \cite{CM:a}). When $|\lambda|=1$,
$U_\lambda$ is unitary and it is easy to see that
\[
T_{\lambda,\vp}=U_{\lambda}T_{\vp_{\ol\lambda,+}},
\]
where $T_{\vp_{\ol\lambda,+}}$ is the Toeplitz operator with symbol
\[
\vp_{\ol\lambda,+}\sim\sum_{n=-\infty}^\infty b_ne^{in\theta},\
b_n={\ol\lambda}^na_n\ \mbox{if}\ n\geq0\ \mbox{and}\ b_n=a_n\
\mbox{if}\ n<0.
\]
The above identity thus puts $T_{\lambda,\vp}$ into the category of
Toeplitz-composition operators when $|\lambda|=1$. When $|\lambda|<1$, $T_{\lambda,\vp}$ also falls into the same category (or the category of {\em weighted composition operators}, see, for example, \cite{CA:a}, \cite{CA:b}, \cite{CK:a}, \cite{GG:a}, \cite{GG:b}, \cite{MZ:a} and \cite{SH:a}) if $\vp$ (or $\ol\vp$) is in $H^2$. On the other hand, however, there is
no reason for us to treat the $\lambda$-Toeplitz operators as a subclass of the
Toeplitz-composition operators since in general, when $|\lambda|<1$, we have
\[
T_{\lambda,\vp}=W_{\vp_+,\tau}+W_{\ol{\vp_-},\ol\tau}^*,
\]
where $W_{\vp_+,\tau}$ is the weighted composition operator
\[
W_{\vp_+,\tau}f:=\vp_+\cdot(f\circ\tau)
\]
and $W_{\ol{\vp_-},\ol\tau}^*$ is the adjoint of $W_{\ol{\vp_-},\ol\tau}$ (For the definition of the adjoint of an operator, see \cite{JC:a}):
\[
W_{\ol{\vp_-},\ol\tau}^*f:=(P(\vp_-f))\circ\tau\ \ \ (W_{\ol{\vp_-},\ol\tau}f=\ol{\vp_-}\cdot(f\circ\ol\tau)),
\]
with $\vp_+=P\vp$, ${\vp_-}=(I-P)\vp$ and $\tau(z)=\lambda z$, $\ol\tau(z)=\ol\lambda z$, $|z|<1$ (Here $I$ is the identity map on $L^2$).

Perhaps the question an operator theorist is most likely to ask about this
subject is that to what degree are $\lambda$-Toeplitz operators and
Toeplitz operators related through the classical results of Toeplitz operators. While it may not
be surprising to know that $\lambda$-Toeplitz operators and Toeplitz
operators can be more or less connected through the classical
Toeplitz operator theory if $|\lambda|=1$ (e.g., Proposition \ref{essential} and Theorem
\ref{spectrum}), they are very different when $|\lambda|<1$. For
instance, one obvious difference between $\lambda$-Toeplitz
operators and Toeplitz operators is that a nontrivial
$\lambda$-Toeplitz operator may be compact, while a Toeplitz
operator is compact if and only if it is the zero operator. In fact,
one can easily show that $T_{\lambda,\vp}$ is compact (or, even
better, in the trace class) if and only if $|\lambda|<1$ or
$\vp\equiv0$, and with some more effort, that $T_{\lambda,\vp}$ is of finite
rank if and only if $\lambda=0$ or $\vp\equiv0$ (See \cite{Ho:h}). It is also worth
mentioning here that if, in addition, $\vp$ is analytic, then some of the
statements about the compactness and the finite rank criteria for
$\lambda$-Toeplitz operators follow directly from Gunatillake's work
in weighted composition operators on Hardy spaces (see, e.g.,
Theorem 1 in \cite{GG:a} and Theorem 2 in \cite{GG:b}).

In this paper, we shall concentrate on investigating the essential spectrum of $T_{\lambda,\vp}$ when $|\lambda|=1$, and we will apply the results to show that the spectrum of $T_{\lambda,\vp}$ equals
\[
\lt\{\mu\in\Bbb C:\mu^q\in\mbox{cl}(\hat\phi_{\ol\lambda+}(\Bbb D))=\mbox{closure of $\hat\phi_{\ol\lambda+}(\Bbb D)$}\rt\}
\]
if $\vp$ (or, what is the same, $\vp_{\ol\lambda,+}$) is analytic and $C^1$, and $\lambda=e^{2i\pi(p/q)}$ with the rational number $p/q$ in lowest terms, where
\[
\phi_{\ol\lambda+}=\prod_{j=0}^{q-1}\vp_{\ol\lambda,+}\circ\ol\tau^j,\
\ \ol\tau(e^{i\theta})={\ol\lambda}e^{i\theta}
\]
and $\hat\phi_{\ol\lambda+}$ is the analytic function, called the {\em Gelfand transform} of $\phi_{\ol\lambda+}$, defined by (See Chap. 6, \cite{RD:a})
\[
\hat\phi_{\ol\lambda+}(z):={1\over2\pi}\int_0^{2\pi}{\phi_{\ol\lambda+}(e^{i\theta})\over1-ze^{-i\theta}}d\theta,\ |z|<1.
\]
This result generalizes, under this $C^1$ restriction for the symbols, a well-known result for the
spectra of analytic Toeplitz operators due to Wintner (See Theorem 7.21, \cite{RD:a} and \cite{AW:a}). As a consequence, we obtain the spectrum for the weighted composition operator $W_{\vp,\rho}$ with $\vp$ being continuously differentiable on $\Bbb T$ and analytic, and $\rho$ being an elliptic analytic automorphism on $\Bbb D$ of finite order (Theorem \ref{application}).

\section{Essential spectrum of $T_{\lambda,\vp}$, $|\lambda|=1$}

We begin with a brief introduction about the Fredholm operators on Hilbert spaces. Let $A$ be a bounded operator on $\cal H$ and $A^*$ be its adjoint. We say that $A$ is {\em Fredholm} (or $A\in\Phi$) if $\dim(\ker A),\ \dim(\ker A^*)<\infty$ and both $A$ and $A^*$ have closed ranges. For $A\in\Phi$, the {\em index} of $A$, denoted by $\mbox{ind}(A)$, is defined by the integer
\[
\dim(\ker A)-\dim(\ker A^*).
 \]
The essential spectrum for a bounded operator $A$, on the other hand, is defined by
\[
\sigma_e(A)=\{\alpha\in\Bbb C:A-\alpha\not\in\Phi\},
\]
and it is a standard property that $\sigma_e(A)=\sigma_e(A+K)$, whenever $K$ is a compact operator (For the Fredholm theory of bounded operators, we refer the readers to \cite{JC:a}).

Now Let $|\lambda|=1$. The following lemma gives the connection between $\sigma_e(T_{\lambda,\vp})$ and $\sigma_e(T_{\phi_{\ol\lambda+}})$ when $\lambda$ is of finite order:
\begin{lma}
For any bounded $\lambda$-Toeplitz operator $T_{\lambda,\vp}$, we
have
\[
T_{\lambda,\vp}^k=U^k_\lambda
T_{\vp_{\ol\lambda,+}\circ\ol\tau^{k-1}}T_{\vp_{\ol\lambda,+}\circ\ol\tau^{k-2}}\cdots
T_{\vp_{\ol\lambda,+}}\ \
(\ol\tau(e^{i\theta})={\ol\lambda}e^{i\theta})
\]
for $k=1,2\cdots$, and as a consequence, we have
\[
\sigma_e(T_{\lambda,\vp})^q=\sigma_e(T_{\lambda,\vp}^q)=\sigma_e(T_{\phi_{\ol\lambda+}})
\]
if $\lambda=e^{2i\pi(p/q)}$ and $\vp_{\ol\lambda,+}\in C(\Bbb T)$ .
\label{identity}
\end{lma}
\pf\ For any $\psi\in L^\infty(\Bbb T)$, we have
\[
\psi U_\lambda f=U_\lambda(({\psi\circ{\ol\tau}})f)\ \ (f\in H^2).
\]
Since obviously $U_\lambda$ commutes with $P$, we have $T_\psi
U_\lambda=U_\lambda T_{\psi\circ\tau}$, and therefore the first
identity follows by  repeated application of this to
\[
T_{\lambda,\vp}^k=\underbrace{U_{\lambda}T_{\vp_{\ol\lambda,+}}U_{\lambda}T_{\vp_{\ol\lambda,+}}\cdots
U_{\lambda}T_{\vp_{\ol\lambda,+}}}_{k}.
\]

Now suppose that $\lambda=e^{2i\pi(p/q)}$ and $\vp_{\ol\lambda,+}\in
C(\Bbb T)$. Since $U_\lambda^q=I$ and $\vp_{\ol\lambda,+}\in C(\Bbb
T)$, we have
\[
T_{\lambda,\vp}^q=T_{\vp_{\ol\lambda,+}\circ\ol\tau^{k-1}}T_{\vp_{\ol\lambda,+}\circ\ol\tau^{k-2}}\cdots
T_{\vp_{\ol\lambda,+}}=T_{\phi_{\ol\lambda+}}+K\
\]
with some compact operator $K$ since Toeplitz operators with continuous symbols commute modulo compact operators (See Proposition 7.22 in \cite{RD:a}). Therefore by the spectral mapping
theorem the proof is complete.\qed\vspace{.05in}\\
Notice that Lemma \ref{identity} implies, in particular,
$\|T_{\lambda,\vp}^k\|_e=\|\phi_{\ol\lambda+}^{(k)}\|_\infty$,
where
\[
\phi_{\ol\lambda+}^{(k)}=\prod_{j=0}^{k-1}\vp_{\ol\lambda,+}\circ\ol\tau^j.
\]
Hence, in terms of ergodic theory, we obtain, through the spectral
radius formula, that (See \cite{PW:a})
\begin{cor}
The essential spectral radius of $T_{\lambda,\vp}$ equals
\[
\sup_{\mu\in M_{\ol\tau}(\Bbb T)}\exp\lt(h_{\ol\tau(}\nu)+\int_{\Bbb
T}\log|\vp_{\ol\lambda,+}|d\nu\rt)
\]
if $\vp_{\ol\lambda,+}\in C(\Bbb T)$, where $M_{\ol\tau}(\Bbb T)$ is the
set of $\ol\tau$-invariant Borel probability measures on $\Bbb T$ and
$h_{\ol\tau}(\mu)$ is the entropy of $\mu$ with respect to ${\ol\tau}$.
\label{essspecradius}
\end{cor}

One important fact pointed out by Lemma \ref{identity} is that if
$\mu\in\sigma_e(T_{\lambda,\vp})$ and $\vp_{\ol\lambda,+}$ is continuous, then $\phi_{\ol\lambda+}-\mu^q$
is not invertible. We will show that
the converse is true if, in addition, $\vp_{\ol\lambda,+}$ is $C^1$ (Proposition \ref{essential}). For this, we need some detailed information from the Fredholm theory for the {\em
Toeplitz-Composition algebra}, i.e., the $C^*$-algebra generated by
the Toeplitz algebra {\bf A} ($C^*$-algebra generated by Toeplitz operators with continuous symbols) and a single composition operator on $H^2$.

Let ${\cal K}$ be the ideal of compact operators on $H^2$. In the 2007 paper \cite{MJ:a}, Michael T. Jury showed that there is an exact sequence
\[
0\longrightarrow{\cal K}\longrightarrow{\bf A}C_\rho\longrightarrow
C(\Bbb T)\rtimes_\rho\Bbb Z\longrightarrow0\,
\]
where ${\bf A}C_\rho=C^*({\bf A},C_\rho)$ is the Toeplitz-Composition algebra generated by {\bf A} and the composition operator $C_\rho$ is induced by an analytic automorphism $\rho$ on $\Bbb D$, and
$C(\Bbb T)\rtimes_\rho\Bbb Z$ is the crossed product of the action
of symbols on $\Bbb T$ (Theorem 2.1, \cite{MJ:a}).
Note that Jury's result generalizes the famous Coburn's exact
sequence (See \cite{LC:a}, \cite{LC:b})
\[
0\longrightarrow{\cal K}\longrightarrow{\bf A}\longrightarrow C(\Bbb
T)\longrightarrow0.
\]
Then he proceeded to show (Theorem 3.1, \cite{MJ:a}) that if $\rho$
is an elliptic automorphism of order $q$, then $T=\sum_{j=0}^{q-1}T_{f_j}C_\rho^j$, with $f_j\in
C^1(\Bbb T)$ for each $j$, is Fredholm if and only
if the $C(\Bbb T)$-valued determinant
\[
h_T=\lt|\begin{array}{cccc}
f_0&f_1&\cdots&f_{q-1}\\
f_{q-1}\circ\rho&f_0\circ\rho&\cdots&f_{q-2}\circ\rho\\
\cdots&\cdots&\cdots&\cdots\\
f_1\circ\rho^{q-1}&\cdots&\cdots&f_0\circ\rho^{q-1}
\end{array}\rt|
\]
is nonvanishing on $\Bbb T$, and, in this case, the Fredholm index $\mbox{ind}(T)$ equals
\[
\mbox{wn}(h_T)={-1\over2\pi iq}\int_{\Bbb T}{dh_T\over h_T},
\]
where $\mbox{wn}(h_T)$ is the {\em winding number} for the curve $h_T(\Bbb T)$.

We now apply Jury's result to get
\begin{prop}
Let $\lambda=e^{2i\pi(p/q)}$ with the rational number $p/q$ in
lowest terms. Suppose, in addition, that $\vp_{\ol\lambda,+}$ is
continuously differentiable. Then the essential spectrum of $T_{\lambda,\vp}$ is
\[
\lt\{\mu\in\Bbb C:\phi_{\ol\lambda+}-\mu^q\ \mbox{is not
invertible}\rt\},
\]
i.e., $\mu\in\sigma_e(T_{\lambda,\vp})$ if and only if $\mu^q\in
\phi_{\ol\lambda+}(\Bbb T)$, where
\[
\phi_{\ol\lambda+}=\prod_{j=0}^{q-1}\vp_{\ol\lambda,+}\circ\ol\tau^j,\
\ \ol\tau(e^{i\theta})={\ol\lambda}e^{i\theta}.
\]
Moreover, we have
$\mbox{ind}(T_{\lambda,\vp}-\mu)=-q^{-1}\mbox{wn}(\phi_{\ol\lambda+}-\mu^q)$.
\label{essential}
\end{prop}
\pf\ Since the unitary $U_{\ol\lambda}$ is in fact
$C_{\ol\tau}$, and
$U_{\ol\lambda}T_{\lambda,\vp}=T_{\vp_{\ol\lambda,+}}$, we see that
$T_{\lambda,\vp}-\mu$ is Fredholm if and only if
$T_{\vp_{\ol\lambda,+}}-\mu C_{\ol\tau}$ is Fredholm. Therefore, by
applying Jury's criterion to $f_0=\vp_{\ol\lambda,+}$, $f_j=0$ for
$j=1,2,\cdots,q-2$, $f_{q-1}=-\mu$, $\rho=\ol\tau$ and $T=T_{\vp_{\ol\lambda,+}}-\mu C_{\ol\tau}$, one has
\begin{eqnarray*}
h_T&=&\lt|\begin{array}{ccccc}
\vp_{\ol\lambda,+}&0&\cdots&\cdots&-\mu\\
-\mu&\vp_{\ol\lambda,+}\circ{\ol\tau}&0&\cdots&0\\
0&-\mu&\cdots&\cdots&0\\
0&\cdots&\cdots&\cdots&0\\
0&0&\cdots&-\mu&\vp_{\ol\lambda,+}\circ{\ol\tau}^{q-1}
\end{array}\rt|\\
&=&\vp_{\ol\lambda,+}\cdot(\vp_{\ol\lambda,+}\circ{\ol\tau})\cdots(\vp_{\ol\lambda,+}\circ{\ol\tau}^{q-1})
+(-1)^{q+1}(-\mu)^q\\
&=&\phi_{\ol\lambda+}-\mu^q,
\end{eqnarray*}
and from this and the Jury's index formula, the index formula for $T_{\lambda,\vp}-\mu$ is derived immediately.\qed\vspace{.05in}\\
\section{Applications to some weighted composition operators}
Let $f\in H^2$ and recall the Gelfand transform of $f$:
\[
\hat f(z):={1\over2\pi}\int_0^{2\pi}{f(e^{i\theta})\over1-ze^{-i\theta}}d\theta,\ |z|<1.
\]
In 1929,  A. Wintner showed in \cite{AW:a} that if $\vp\in L^\infty$ is analytic and $\hat\vp$ is bounded on $\Bbb D$ (or, simply, $\vp\in H^\infty$), then the spectrum of the Toeplitz operator $T_{\vp}$ is $\mbox{cl}(\hat\vp(\Bbb D))$. Notice that $T_\vp=T_{1,\vp}$ and the order of 1 is one. We now generalize this result to $\lambda$-Toeplitz operators, with an additional assumption that $\vp$ is smooth:
\begin{thm} Let $\lambda=e^{2i\pi(p/q)}$ with the rational $p/q$ in
lowest terms. Then $\sigma(T_{\lambda,\vp})$, the spectrum of
$T_{\lambda,\vp}$, equals
\[
\lt\{\mu\in\Bbb C:\mu^q\in\mbox{cl}(\hat\phi_{\ol\lambda+}(\Bbb D))\rt\}
\]
if $\vp$ is analytic and continuously differentiable.
\label{spectrum}
\end{thm}
\pf\ We will prove the theorem by showing that $T_{\lambda,\vp}-\mu$ is invertible if and only if $T_{\phi_{\ol\lambda+}}-\mu^q$ is. The conclusion then follows directly from Wintner's theorem for analytic Toeplitz operators.

It is well-known in the Fredholm theory of Toeplitz operators that if $T_\psi$ is Fredholm, then $T_\psi$ is invertible if and only if $\mbox{ind}(T_\psi)=0$, and $\mbox{ind}(T_\psi)=-\mbox{wn}(\psi)$ if $\psi$ is continuous (See Corollary 7.25, \cite{RD:a}).
On the other hand, by elementary operator theory, we know that if $A$ is not invertible, then only one of the following three possibilities may occur: (1) $A\not\in\Phi$ (2) $A\in\Phi$, $\mbox{ind}(A)\not=0$ (3) $A\in\Phi$, $\mbox{ind}(A)=0$, but not invertible. Hence the spectrum of $T_{\phi_{\ol\lambda+}}$ equals
\begin{eqnarray*}
&&\lt\{\alpha:T_{\phi_{\ol\lambda+}-\alpha}\not\in\Phi\rt\}\cup\lt\{\alpha:T_{\phi_{\ol\lambda+}-\alpha}\in\Phi,\ \mbox{ind}(T_{\phi_{\ol\lambda+}-\alpha})\not=0\rt\}\\
&=&\sigma_e(T_{\phi_{\ol\lambda+}-\alpha})\cup\lt\{\alpha:T_{\phi_{\ol\lambda+}-\alpha}\in\Phi,\ \mbox{wn}(\phi_{\ol\lambda+}-\alpha)\not=0\rt\}
\end{eqnarray*}
since $\phi_{\ol\lambda+}$ is continuous (Note that $T_\psi-\alpha=T_{\psi-\alpha}$).
So, by the Fredholm criterion and the index formula in Proposition \ref{essential},
the only thing we need to show is that if $T_{\lambda,\vp}-\mu\in\Phi$ and
$\mbox{ind}(T_{\lambda,\vp}-\mu)=0$, then $T_{\lambda,\vp}-\mu$ is
invertible.

Since $\vp_{\ol\lambda,+}=\vp\circ\ol\tau$ if $\vp$ is analytic, $\vp\in H^\infty$ if and only if
$\vp_{\ol\lambda,+}\in H^\infty$. Therefore,
\[
T_{\lambda,\vp}^q=T_{\vp_{\ol\lambda,+}\circ\ol\tau^{k-1}}T_{\vp_{\ol\lambda,+}\circ\ol\tau^{k-2}}\cdots
T_{\vp_{\ol\lambda,+}}=T_{\phi_{\ol\lambda+}}\ \mod{\cal K}
\]
in Lemma \ref{identity} is actually an equality.
Thus we have
\[
\sigma(T_{\lambda,\vp})^q=\sigma(T_{\phi_{\ol\lambda+}}).
\]
Now by Proposition \ref{essential}, if $T_{\lambda,\vp}-\mu\in\Phi$
and $\mbox{ind}(T_{\lambda,\vp}-\mu)=0$, then
$\phi_{\ol\lambda+}-\mu^q$ is invertible and
$\mbox{wn}(\phi_{\ol\lambda+}-\mu^q)=0$. This means that,
$T_{\phi_{\ol\lambda+}}-\mu^q$ must be invertible, which implies
immediately that $\mu\not\in\sigma(T_{\lambda,\vp})$.\qed\vspace{.05in}\\
{\bf Remark.} It is not difficult to see the connection between the set $\hat\phi_{\ol\lambda+}(\Bbb D)$ and $\phi_{\ol\lambda+}$ since for $f\sim\sum_0^\infty a_ne^{i\theta}\in H^2$, the Taylor series of $\hat f$ at the origin is $\sum_0^\infty a_n z^n$, $|z|<1$. Furthermore, by a theorem of Fatou, $\hat f(re^{i\theta})\ra f(e^{i\theta})$ a.e. $\theta$ and in $L^2$ (See \cite{RD:a} and \cite{KH:a}). Hence we may consider $f$ as the ``boundary value" of $\hat f$. Let us use this understanding in the following example:

Let $r(A)$ and $r_e(A)$ denote the spectral radius and essential spectral radius of $A$, respectively. Let $\lambda=e^{i2\pi/3}$, and
consider $\vp(e^{i\theta})=e^{i\theta}-\sqrt[3]{2}$. Then
$\vp_{\ol\lambda,+}(e^{i\theta})=e^{i4\pi/3}e^{i\theta}-\sqrt[3]{2}$ and
$\phi_{\ol\lambda+}(e^{i\theta})=e^{i3\theta}-2$. Therefore
\[
r_e(T_{\lambda,\vp})=r(T_{\lambda,\vp})=\sqrt[3]{3}<1+\sqrt[3]{2}=\|\vp_{\ol\lambda,+}\|_\infty=\|T_{\lambda,\vp}\|.
\]
This means that the identity
\[
r(T_\vp)=r_e(T_\vp)=\|T_\vp\|,\ \vp\in C(\Bbb T)
\]
for Toeplitz operators is no longer valid for $\lambda$-Toeplitz
operators.

The same example also shows that a well-known result about the
connectedness of the spectrum of Toeplitz operators (due to Windom)
no longer holds for $\lambda$-Toeplitz operators since $\hat\phi_{\ol\lambda+}(z)=z^3-2$ and so by Theorem
\ref{spectrum},
\[
\sigma(T_{\lambda,\vp})=\{\mu\in\Bbb
C:\mu^3\in\mbox{cl}(\hat\phi_{\ol\lambda+}(\Bbb D))\}=\{\mu\in\Bbb
C:|\mu^3+2|\leq1\}.
\]
Clearly, this set is not connected.

We finish with an application of Theorem \ref{spectrum} to certain weighted composition operators:
\begin{thm}
Let $\vp\in H^\infty\cap C^1$ and $\rho$ be an elliptic analytic automorphism of $\Bbb D$ of order $q$. Then the spectrum of the weighted composition $W_{\vp,\rho}$ is
\[
\lt\{\mu\in\Bbb C:\mu^q\in\mbox{cl}(\hat\psi(\Bbb D))\rt\},
\]
where
\[
\psi=\prod_{k=0}^{q-1}(\vp\circ\zeta^{-1})\circ\ol\tau^k,\
\]
$\zeta$ is an analytic automorphism of $\Bbb D$ such that $\zeta^{-1}(\rho(\zeta(z)))=\lambda z$ for some $\lambda\in\Bbb T$ with order $q$, and $\ol\tau(z)=\ol\lambda z$.\label{application}
\end{thm}
\pf\ Since $\rho$ is an elliptic analytic automorphism of $\Bbb D$ of order $q$, we can find analytic automorphism $\zeta$ of $\Bbb D$ such that $\zeta^{-1}(\rho(\zeta(z)))=\lambda z$, $|z|<1$ for some $\lambda\in\Bbb T$ with order $q$. Now since the composition operator $C_\zeta$ is invertible and $C_\zeta^{-1}=C_{\zeta^{-1}}$, we have
\[
C_\zeta^{-1}W_{\vp,\rho}C_\zeta=C_{\zeta^{-1}}W_{\vp,\rho}C_\zeta=W_{\vp\circ\zeta^{-1},\tau},\ \tau(z)=\lambda z,
\]
which means that $W_{\vp,\rho}$ is similar to $W_{\vp\circ\zeta^{-1},\tau}$.

Now since $\vp\in H^\infty\cap C^1$ implies $\vp\circ\zeta^{-1}\in H^\infty\cap C^1$, $W_{\vp\circ\zeta^{-1},\tau}$ is the same with $T_{\lambda,\vp\circ\zeta^{-1}}$, and hence the result follows from Theorem \ref{spectrum} since we also have $f_{\ol\lambda,+}=f\circ\ol\tau$ if $f$ is analytic.\qed

Chih Hao Chen\\
Department of Applied Mathematics\\
National Sun Yat-Sen University\\
Kaohsiung, Taiwan 80424\\
m972040018@student.nsysu.edu.tw\vspace{.1in}\\
Po Han Chen\\
Department of Applied Mathematics\\
National Sun Yat-Sen University\\
Kaohsiung, Taiwan 80424\\
m972040010@student.nsysu.edu.tw
\vspace{.1in}\\
Mark C. Ho\\
Department of Applied Mathematics\\
National Sun Yat-Sen University\\
Kaohsiung, Taiwan 80424\\
hom@mail.nsysu.edu.tw\vspace{.1in}\\
Meng Syun Syu\\
Department of Applied Mathematics\\
National Sun Yat-Sen University\\
Kaohsiung, Taiwan 80424\\
m972040015@student.nsysu.edu.tw

\end{document}